\documentclass[12pt]{article}
\usepackage{amsmath}
\usepackage{amssymb}

\usepackage{graphicx}

\let\rho=\varrho

\def\qed{\hfill$\Box$\par\medskip\par\relax}

\def\1#1{{\bf 1}{\{#1\}}}
\newcommand{\eps}{\varepsilon}
\newcommand{\Z}{{\mathbb Z}}
\newcommand{\N}{{\mathbb N}}
\newcommand{\R}{{\mathbb R}}
\newcommand{\C}{{\mathbb C}}
\newcommand{\IP}{{\mathbb P}}
\newcommand{\E}{{\mathbb E}}
\newcommand{\V}{{\mathcal V}}
\newcommand{\W}{{\mathcal W}}

\newcommand{\G}{{\mathcal G}}
\newcommand{\Perc}{{\tt Perc}}
\newcommand{\calP}{{\mathcal P}}

\newcommand{\calE}{{\mathcal E}}
\newcommand{\calI}{{\mathcal I}}

\newcommand{\Q}{{\mathbb Q}}
\let\phi=\varphi

\newcommand{\calA}{{\mathcal A}}

\newcommand{\pc}{{\vec p_c}}

\newcommand{\Si}{{\Sigma}}
\newcommand{\be}{{\beta}}
\newcommand{\h}{{\eta}}
\newcommand{\al}{{\alpha}}
\newcommand{\la}{{\lambda}}
\newcommand{\lh}{{\hat \lambda}}
\newcommand{\8}{{\infty}}
\newcommand{\Card}{\mathop{\rm Card}\nolimits}

\renewcommand{\Re}{\mathop{\rm Re}\nolimits}
\renewcommand{\Im}{\mathop{\rm Im}\nolimits}

\allowdisplaybreaks

\newtheorem{theo}{Theorem}[section]
\newtheorem{lmm}[theo]{Lemma}

\newtheorem{prop}[theo]{Proposition}

\newtheorem{rem}[theo]{Remark}
\newtheorem{conj}[theo]{Conjecture}

\title{The number of open paths in an oriented $\rho$-percolation model}
\author{Francis~Comets\thanks{Partially supported by CNRS (UMR 7599
``Probabilit{\'e}s et Mod{\`e}les Al{\'e}atoires'') and 
ANR Polintbio}$^{~,1}$ \and
 Serguei~Popov\thanks{Partially supported by CNPq (302981/02--0),
 and by the ``Rede Mate\-m\'atica Brasil-Fran\c{c}a''}$^{~,2}$
\and Marina Vachkovskaia\thanks{Partially supported 
         by CNPq (304561/2006--1  and 200460/2006--4)}$^{~,3}$}

\begin{document}

\maketitle

{\footnotesize
\noindent $^{~1}$Universit{\'e} Paris 7, UFR de Math{\'e}matiques,
case 7012, 2, place Jussieu, F--75251 Paris Cedex 05, France\\
\noindent e-mail: \texttt{comets@math.jussieu.fr},
\noindent url: \texttt{http://www.proba.jussieu.fr/$\sim$comets}

\smallskip
\noindent $^{~2}$Instituto de Matem{\'a}tica e Estat{\'\i}stica,
Universidade de S{\~a}o Paulo, rua do Mat{\~a}o 1010, CEP 05508--090,
S{\~a}o Paulo SP, Brasil\\
\noindent e-mail: \texttt{popov@ime.usp.br}, 
\noindent url: \texttt{http://www.ime.usp.br/$\sim$popov}

\smallskip
\noindent $^{~3}$Departamento de Estat\'\i stica, Instituto de Matem\'atica,
Estat\'\i stica e Computa\c{c}\~ao Cien\-t\'\i\-{}fica,
Universidade de Campinas,
Caixa Postal 6065, CEP 13083--970, Campinas SP, Brasil\\
\noindent e-mail: \texttt{marinav@ime.unicamp.br},
\noindent url: \texttt{http://www.ime.unicamp.br/$\sim$marinav}
}

\begin{abstract}
We study the asymptotic properties of the
number of open paths of length~$n$ in an oriented $\rho$-percolation model.
We show that this number is $e^{n\alpha(\rho)(1+o(1))}$ as $n \to \8$.
The exponent $\alpha$ is deterministic, it can be expressed in terms
of the free energy of a polymer model, 
and it can be explicitely computed in some range of the parameters.
Moreover, in a restricted  range of the parameters, we even show
that the number of such paths is $n^{-1/2} W e^{n\alpha(\rho)}(1+o(1))$
for some nondegenerate random variable $W$.
We build on connections with the model of  directed
polymers in random environment, and we use  techniques and results
developed in this context.
%% Our approach is based on several techniques developed for directed
%% polymers in random environment. 
\\[.3cm]{\bf Short Title:} Number of open paths in oriented $\rho$-percolation
\\[.3cm]{\bf Keywords:} Oriented percolation, $\rho$-percolation,
 directed polymers in random environment.
\\[.3cm]{\bf AMS 2000 subject classifications:} 60K37
\end{abstract}

% {\bf TO DO LIST}
% \begin{enumerate}
% \item faire $\rho \nearrow 1$ ?
% %% $ http://en.wikipedia.org/wiki/Directed_percolation$
% %\item Abstract
% \item Intro
% \end{enumerate}

\section{Introduction and results}
\label{s_intro}

\subsection{Introduction}
In this paper we study the number of open paths in an oriented
$\rho$-percolation model in dimension $1+d$, or, equivalenly, the number 
of $\rho$-open
path in an oriented percolation model. 
Consider the graph $\Z_+\times\Z^d$, with $\Z_+=\{0,1,2,3,\ldots\}$, 
and fix some parameter~$p\in(0,1)$. 
To each site of this graph except the origin, assign a variable taking 
value~$1$ with probability~$p$ and~$0$
with probability $1-p$, independently of the other sites. 
An oriented (sometimes also called \emph{semi-oriented}) path
of length~$n$ is a sequence $(0,x_0), (1,x_1), (2,x_2),\ldots,(n,x_n)$,
where $x_0=0$ and $x_i,x_{i+1}$ are neighbours in~$\Z^d$, $i=0,\ldots,n-1$:
viewing  the first coordinate as time, 
one can think of such path as a path of the $d$-dimensional simple random walk.
Fix another parameter $\rho\in[0,1]$; the concept of $\rho$-percolation
was introduced by Menshikov and Zuev in \cite{MZ}, as the occurence
of an infinite length path with asymptotic density of 1s larger of equal to 
$\rho$. As in classical percolation, this probability of this 
event is subject to a dychotomy \cite{MZ}
according to $p$ larger or smaller than some critical threshold, which
was later studied by Kesten and Su \cite{KS} in the asymptotics of 
large dimension.

In the present paper, we discuss paths of finite length $n$, in the limit
$n \to \8$. An oriented path of length~$n$
is called $\rho$-open, if the proportion of $1$s in it is at least~$\rho$.
From standard percolation theory it is known 
%%It is known 
%%(see e.g.~\cite{KS,MZ}) 
that for large~$p$ there are 
$1$-open oriented paths with 
%uniformly (in $n$) positive 
nonvanishing
probability, and 
from ~\cite{MZ}
that for any~$p$ one can 
find 
%% a small enough~$\rho$ 
$\rho$ larger than $p$
such that, almost surely, 
%%, with probability that converges to~$1$ as $n\to\infty$, 
there are $\rho$-open oriented paths
for large $n$. However, the question of how many such
paths of length~$n$ can be found in a typical situation, was still unaddressed 
in the literature. When 
%% preparing 
finishing this manuscript, we have learned
of the related work~\cite{KSid}. 

In this paper, we prove that the number of different $\rho$-open
paths of length~$n$ behaves like $e^{n\alpha(\rho)(1+o(1))}$,
where the exponent $\alpha(\rho)$ is deterministic and, of course,
also depends on~$p$ and~$d$. We prove that the function  $\alpha(\cdot)$ 
is the negative convex conjugate of the free energy of
{\it directed polymers in random environment}. This model has attracted
a lot of interest in recent years, leading to a better -- although very 
incomplete -- understanding. We will extensively use the current knowledge of 
thermodynamics of the polymer model, and 
the reader is refered to \cite{CSYrev} for a recent survey. 
This will allow us to obtain, when~$d \geq 3$, the explicit expression 
for $\alpha(\rho)$
in a certain range of values for $\rho$ depending 
on the parameters~$p$ and~$d$. The reason for this remarkable fact is the
existence of the so-called ``weak-disorder region'' in the polymer model,
discovered in \cite{IS} and \cite{Bolt}:
this reflects here into a parameter region where the number of paths is
of the same order as its expected value. 
\medskip

At this point the reader may be tempted to use first and second moment methods
 to estimate the number of paths. 
The first moment is easily computed, and serves as an upper bound in complete 
generality. 
The second moment is more difficult to analyse. However, it can be checked 
that in large dimension and for density close to the parameter $p$ of 
the Bernoulli, the ratio second-to-first-squared remains bounded 
in the limit of an infinitely long path. This means that, under these 
circumstances, the upper bound gives the right 
order of magnitude with a positive probability. However, 
(i) this method does not tell us anything on $\alpha$ for general parameters, 
(ii) it fails to keep track of the correlation between counts for 
different values of the density.

Our strategy will be quite different. We will study
the {\it moment generating function} of the number of paths, which is not
surprising in such a combinatorial problem. The point is that 
the  moment generating function is simply 
the partition function of the directed polymer in random environment.
This is a well-known object in statistical physics,
its logarithmic asymptotics is well studied, and is given by the free 
energy. From the existence and known properties of the free energy, 
we will derive the existence of $\alpha$ and its expression in 
thermodynamics terms. 
In the course of our analysis we will prove that the free energy,
a convex function of the inverse temperature, is in fact strictly convex. 
This property is new and interesting for the polymer model.

Moreover, in a more restricted  range of  values for $\rho$, we even obtain
an equivalent for the number of paths which achieves exactly a given density
of 1s. This is clearly a very sharp estimate, that we obtain by 
using the power of complex 
analysis, and convergence of the renormalized moment generating function
in the sense of analytic functions.
Certainly a naive moments method cannot lead to such an equivalent.

\subsection{Notations and results}
Now, let us define the model formally.
Let $\h(t,x), t = 1,2,\ldots, x \in \Z^d$ be a sequence of independent 
identically distributed Bernoulli random variables, with common parameter
$p \in (0,1)$, $\IP(\h(t,x)=1)=p=1-\IP(\h(t,x)=0)$.  
We denote by
$(\Omega, \calA, \IP)$ the probability space where this sequence is defined.
The vertex $(t,x)$
is {\it open} if $\h(t,x)=1$ and {\it closed} in the opposite case.
A nearest neighbour 
path~$S$ in $\Z^d$ of length $n$ ($1 \leq n \leq \8$) is a sequence
$S=(S_t; t=0,\ldots, n), S_t\in \Z^d, S_0=0, \|S_t-S_{t-1}\|_1
=1$ for $t=1,\ldots, n$. We denote by $\calP_n$ the set of such paths $S$, 
and by $\calP_\8$ 
the set of infinite length nearest neighbour paths. For $S\in\calP_n$, let
\begin{equation} \label{eq:energy}
H_n(S)=\sum_{t=1}^n \h(t,S_t)
\end{equation}
be the number of open vertices along the path $S$. 

In oriented percolation, one is concerned with the event that there exists 
an infinite open path $S$, i.e. 
$$
\Perc=\big\{ \text{there exists }S \in \calP_\8: \h(t,S_t)=1  \text{ for all }t \geq 1 \big\}\;.
$$ 
It is well known \cite{Dur,Gr} that there exists $\pc(d) \in (0,1)$, called the  
critical percolation threshold,  such that
 %    $$
 %    \IP( \Perc)
 %    \left\{
 %    \begin{array}{c}
 %      >0 \\ =0
 %    \end{array} 
 %    \right. \quad
 %    {\rm \ if \ } \quad p \;
 %    \left\{
 %    \begin{array}{c}
 %     > \pc(d)\\
 %     < \pc(d)
 %    \end{array} 
 %    \right.
 %    $$
\begin{equation}
  \label{eq:11062}
\IP( \Perc) 
\left\{
\begin{array}{ccc}
  > 0 &
{\rm \ if \ } & p 
 > \pc(d),\\  =0 &
{\rm \ if \ } & p  < \pc(d).
\end{array} 
\right.  
\end{equation}
For $\rho \in (p,1]$, Menshikov and Zuev \cite{MZ} introduced 
$\rho$-percolation as the event that there exists an infinite path $S$
with asymptotic proportion at least $\rho$ of open sites,
$$
\text{$\rho$-}\Perc=\big\{ \text{there exists }S \in \calP_\8: 
\liminf_{n \to \8} H_n(S)/n \geq \rho
 \big\}
\;.$$ 
They showed that there also exists a threshold $\pc(\rho, d)$ such
that (\ref{eq:11062}) holds with $\text{$\rho$-}\Perc$ instead of $\Perc$
(with the probability of $\text{$\rho$-}\Perc$ being equal to 1 when
$p>\pc(\rho, d)$). Very little has been proved for  $\rho$-percolation.
 The asymptotics 
of $\pc(\rho, d)$ for large $d$ are obtained in \cite{KS} at first order,
showing that $d^{1/\rho}\pc(\rho, d)$ has a limit as $d \to \8$, and 
that the limit is different from 
the analogous quantity for $d$-ary trees.
As mentioned in this reference, the equality $\pc(1, d)=\pc( d)$ 
follows from Theorem 5 of \cite{Lee}.

In this paper we are interested in the number of oriented paths of length~$n$ 
which have exactly~$k$ open vertices ($k \in \{0,\ldots, n\}$),
$$
Q_n(k) =\Card 
\big\{ S \in \calP_n: H_n(S)=k \big\}
$$
($\Card A$ denotes the cardinality of $A$) and 
%its 
%% cumulative count 
%tail 
the related quantity
given,
for $\rho \in [0,1]$, by
$$
R_n(\rho) =
\left\{
\begin{array}{cc}
\Card 
\big\{ S \in \calP_n: H_n(S) \geq n \rho \big\}\;, & \rho \geq p,\\
\Card 
\big\{ S \in \calP_n: H_n(S) \leq n \rho \big\}\;,\vphantom{\int\limits^N} & \rho < p.
\end{array} 
\right.
$$
%(observe that, by symmetry, $R_n(\rho)$ has the same law as ...).
Note that $Q_n(k), R_n(\rho)$ are random variables, that
$R_n(\rho)=\sum_{k \geq n\rho}Q_n(k)$ when $\rho \geq p$,
and that 
$\Perc=\bigcap_n \{ Q_n(n)\geq 1\}=\bigcap_n \{R_n(1)\geq 1\}$.

\begin{figure}
\centering
\includegraphics{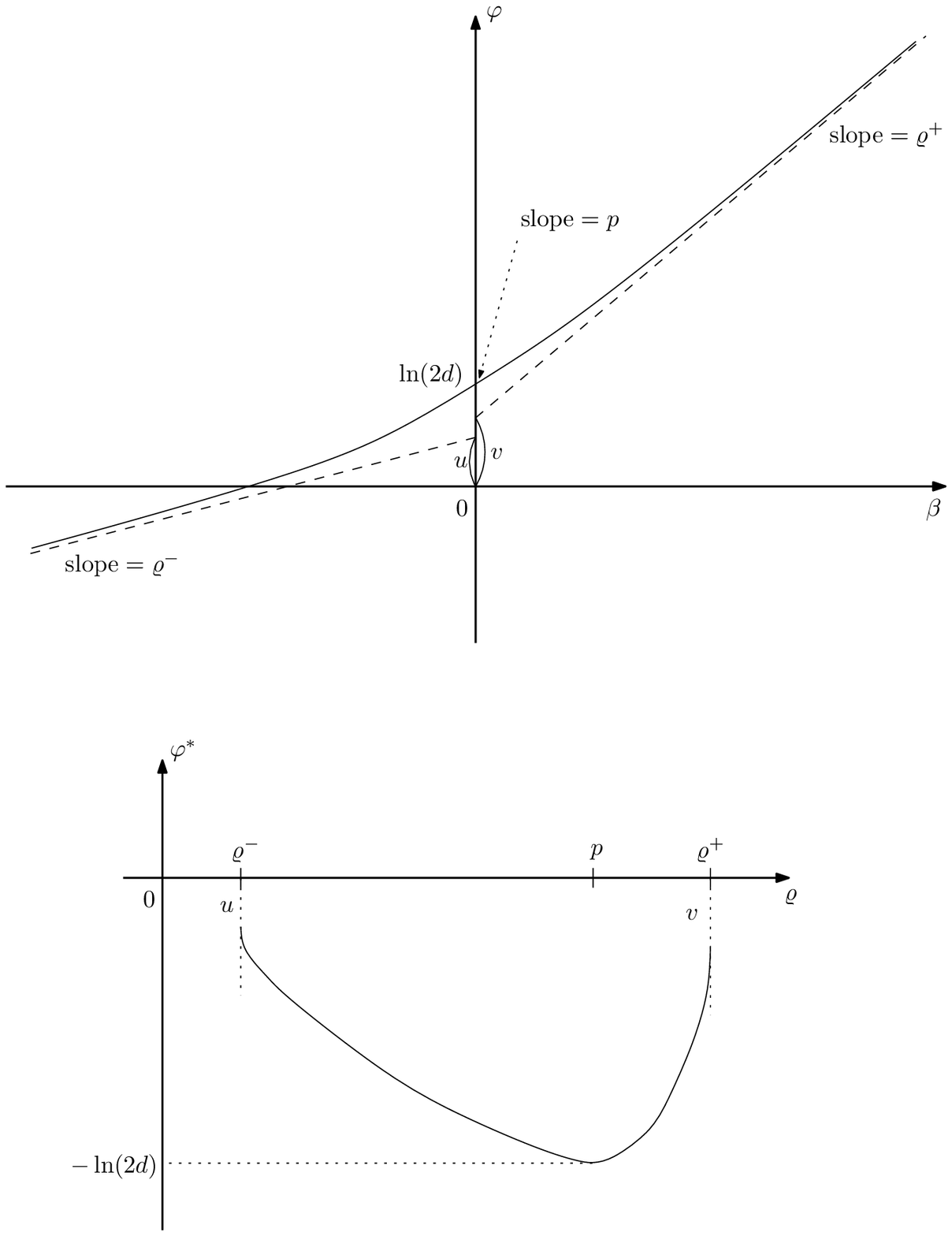}
\caption{The function $\phi$ and its Legendre transform~$\phi^*$}
\label{fig:phi}
\end{figure}

In this paper we relate these quantities to the model of directed polymers
in random environment. Central in this model is the 
(unnormalized) partition function
$Z_n=Z_n(\be, \h)$
at inverse temperature $\be \in \R$
in the environment $\h$ given by
$$
Z_n= \sum_{S \in \calP_n} \exp\{ \be H_n(S) \}\;.
$$
By subadditive arguments one can prove that
\begin{equation} \label{def:phi}
\phi(\be) = \lim_{n \to \8} \frac{1}{n} \E \ln Z_n 
\end{equation}
exists in $\R$ ($ \E$ is the expectation under $\IP$), and by concentration 
arguments, that the event $\Omega_0(\be)$ defined by
\begin{equation} \label{def:phips}
 \Omega_0(\be) = \Big\{\lim_{n \to \8} \frac{1}{n} \ln Z_n = \phi(\be)\Big\}
\end{equation}
has full measure, $\IP(\Omega_0(\be))=1$, see e.g.\ \cite{CSY}.
The function $\phi$ is called the {\it free energy}, it is a
non-decreasing and convex function of  $\be$. Its Legendre 
conjugate
\begin{equation} \label{def:phi*}
\phi^*(\rho) = \sup\{  \be \rho - \phi(\be); \be \in \R\}\;,
\end{equation}
is a convex, lower semi-continuous 
function from $[0,1]$ to $\R \cup \{+\8\}$, such that 
$ \phi^*(\rho)\geq \phi^*(p)= -\ln (2d)$
(indeed,  $\phi'(0)=p$, as it will be shown later). Legendre 
convex duality is better understood by taking a glance at the
graphical construction, e.g.\ figures 2.2.1 and 2.2.2 in \cite{DZ};
here, on Figure~\ref{fig:phi} we illustrate how the functions~$\phi$
and~$\phi^*$ typically look in our situation. 
The existence of the so-called time constants, 
$$
\rho^+=\lim_{n \to \8} \max_{S\in\calP_n} \frac{H_n(S)}{n}\;,\quad
 \rho^-=\lim_{n \to \8} \min_{S\in\calP_n} \frac{H_n(S)}{n}\;, \quad \text{ $\IP$-a.s.}
$$
can be obtained by specifying a direction for the ending point $S_n$, 
which allows using subadditive arguments \cite{Kest},
and then summing over the possible directions. However we give here 
a short proof in the spirit of this work. Since
$
\exp \{ \beta  \max_{S\in\calP_n} {H_n(S)}\}
\leq
Z_n
\leq 
(2d)^n \exp \{ \beta  \max_{S\in\calP_n} {H_n(S)}\}
$
, we have
$$
\frac{1}{n\beta} \ln Z_n-\frac{1}{\beta}\ln(2d) \leq
 \max_{S\in\calP_n} \frac{H_n(S)}{n}
\leq 
\frac{1}{n\beta} \ln Z_n\;.
$$
Taking the limits $n \to \8$ and then $\beta \to +\8$, we see that 
$\rho^+$ is well-defined, and is in fact equal to the slope
$\lim_{\beta \to +\8} \phi(\beta)/\beta$ 
of the asymptotic direction of $\phi$ at $+\8$.
From standard properties of convex duality, the range of the derivative 
$(d/d\be)(1/n)\ln Z_n(\be)$
a.s.\ converges to $ [\rho^-,\rho^+]$  in 
Hausdorff distance,
and $\phi^*(\rho)<+\8$ if and only if $\rho \in [\rho^-,\rho^+]$.
For such $\rho$, we have $\phi^*(\rho) \leq 0$.
\begin{theo} \label{th:Rn}
For all $\rho \in [0,1]$ with 
$ \rho \neq \rho^+, \rho^-,$ the following limit 
\begin{equation} \label{th:eqRn}
\alpha(\rho)= \lim_{n \to \8} \frac{1}{n} \ln R_n(\rho)
\end{equation}
exists $\IP$-a.s.\ (possibly assuming the value $-(\infty)$),  and is given by
$$
\alpha(\rho) = - \phi^*(\rho)\;.
$$
\end{theo}
Clearly, $\al$ is concave, with values in $[0,\ln(2d)]\cup\{-\8\}$ and
$\al(p)=\ln(2d)$.
Note that, for all $\rho \in (\rho^-,\rho^+)$ and almost every
$\h$, 
$$
 R_n(\rho) = \exp n[\al(\rho)+o(1)]\;,\quad \text{as } \; n \to \8 \;.
$$

\goodbreak

\begin{rem}  \label{rem:cont:al}  \end{rem}  \vspace{-3mm}
%%As we will see below (Theorem~\ref{theo:phistrictcvx}), t
By convexity the function $\al$ is continuous on
$(\rho^-,\rho^+)$. For now, it is not clear to us whether 
the limit $\alpha(\rho^+-)$ should be equal to~$0$ in the case 
$p \leq \pc(d)$. In the case $p > \pc(d)$, it is possible to show by 
subadditive arguments that, conditionally on percolation, 
the limit $\alpha(1)$ in (\ref{th:eqRn}) exists and is positive, but
it  is not clear to us  whether $\alpha$ is continuous at 1.

\medskip

Let
$$
\la(\be) = \ln \E e^{\be \h(t,x)} = \ln\left[ 1 + p(e^\be-1)\right]\;,\qquad 
\lh(\be)=\la(\be)+\ln(2d)\;,
$$
then $\E Z_n= \exp \{n \lh(\be)\}$. A direct computation shows that the Legendre 
conjugate $\lh^*(\rho) = \sup\{  \be \rho - \lh(\be); \be \in \R\}$
of $\lh$ is equal to 
\begin{eqnarray} \label{eq:lh*}
\lh^*(\rho) &=& - \ln(2d) + \rho \ln \frac{\rho }{p} + (1-\rho) 
\ln \frac{1-\rho }{1-p}\\ \nonumber
&=&  \rho \ln \frac{\rho }{2dp} + (1-\rho) 
\ln \frac{1-\rho }{2d(1-p)}.
\end{eqnarray}
The  function $(-\lh^*)$ is important for understanding the rate
$\al$. (We recall that both functions depend on $p$, but we don't write
explicitely the dependence.) 
Note that these two functions coincide at $\rho=p$ and take the 
value $\ln (2d)$.
\begin{theo}  \label{th:al-lh} Let $p \in (0,1)$.

  \begin{enumerate}
  \item We have \emph{the annealed bound}: For all $\rho$,
    \begin{equation} \label{eq:annbound}
    \al(\rho) \leq -\lh^*(\rho)\;.  
    \end{equation}
  \item The function $ \al(\rho) + \lh^*(\rho)$ is nonincreasing for
$\rho \in [p, \rho^+)$ and  is nondecreasing for
$\rho \in (\rho^-,p]$.
  \item The set 
\begin{equation} \label{def:Vp}
\V(p)=\{\rho \in (0,1): \al(\rho) = -\lh^*(\rho)\}
\end{equation}
 is an interval containing $p$ (here, ``interval'' is understood in broad 
sense,
i.e., it can reduce to the single point $\{p\}$).
  \item In dimension $d=1$, $\V(p)=\{p\}$, i.e. the inequality in 
    (\ref{eq:annbound}) is strict
    for all $\rho \neq p$.
  \item In dimension $d\geq 3$, $\V(p)$ contains a neighborhood 
    of $p$.
  \item Let $d\geq 3$, and $\pi_d$ be the probability
for the $d$-dimensional simple random walk to ever return to the 
starting point.
When $p>\pi_d$,
    then  $[p,1) \subset \V(p)$, so that the equality holds in 
    (\ref{eq:annbound}) for all $\rho \in
    [p,1)$. Similarly,  when $p<1-\pi_d$,  then $(0,p] \subset \V(p)$, 
so that the equality holds for all $\rho \in
    (0,p]$.
  \item In dimension $d\geq 2$, if $p < (1/2d)$, then $\sup \V(p)<1$.  
Similarly,  if $p > 1-(1/2d)$, we have
  $\inf \V(p)>0$.
  \end{enumerate}
\end{theo}

\begin{rem}  \label{rem:al-lh}  \end{rem}  \vspace{-3mm}
(i) The annealed bound comes from the first-moment method, and  most of the
results
stating that the equality $\al(\rho) = -\lh^*(\rho)$ holds, are derived
from the second-moment method. 

\noindent
(ii) By transience of the  random walk in dimension $d \geq 3$, we have $\pi_d<1$.
In fact, $\pi_3=0.3404\ldots > \pi_4 > \pi_5 \ldots$ \cite[page 103]{Spitz}.
In particular, for $p \in (\pi_d, 1-\pi_d)$, 
we have $\al(\rho) = -\lh^*(\rho)$ for all $\rho\in (0,1)$.

\medskip

\begin{figure}
\centering
\includegraphics{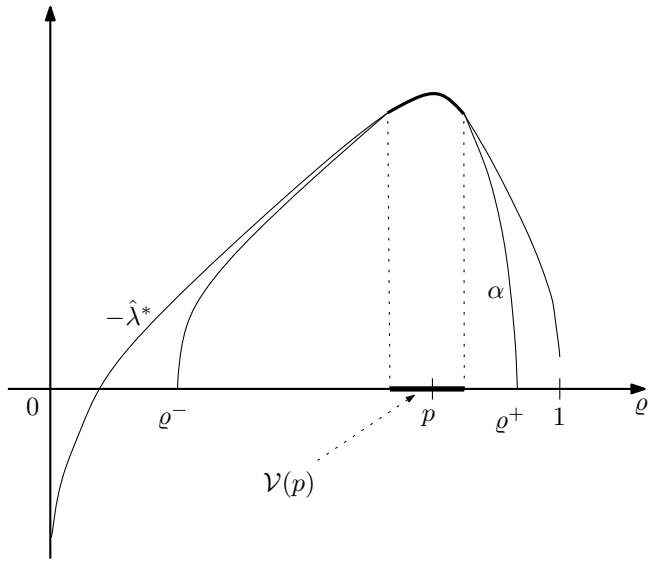}
\caption{Typical behaviour of the function $\alpha$ when $d\geq 3$}
\label{fig:lambda_alpha}
\end{figure}

The following property of the free energy $\phi$ 
of the directed polymer is interesting and seems to be new.

\begin{theo} \label{theo:phistrictcvx}
The function $\phi$ is strictly convex on $\R$, and the functions 
$\phi^*$ and $\al$
are differentiable in the interior of their domains.
\end{theo}

We will obtain much sharper results for large dimension and $\rho$'s not
too far from $p$. The reason is that the partition function $Z_n$ behaves
smoothly as $n \nearrow \8$. The almost-sure limit 
$$
W_\8(\be) = \lim_{n \to \8} Z_n(\be) e^{-n \lh(\be)}
$$
exists for all $\be$, since the sequence is a positive 
 $(\G_n)_n$-martingale, where $\G_n=\sigma\{\h(t,x); 
t \leq n, x \in \Z^d\}$. 
So, let us now concentrate on the case of large dimension, $d \geq 3$. 
When $\be$ belongs to some neighborhood of the origin
(known as the weak disorder region), the limit $W_\8$ is strictly positive
a.s.  In a smaller neighborhood of the 
origin,  the limit can be expressed as a (random) perturbation 
series in ${\mathcal L}^2$
\cite{Sinai}. Moreover, 
the convergence holds in much stronger sense, namely, in the sense of 
analytic functions \cite{CY}. We will use strong tools from complex analysis,
as it is classically done to obtain limit theorems for sums of random 
variables \cite{Petrov}.

\begin{theo} \label{theo:sharp1}
Assume $d \geq 3$. There exist a  neighborhood $U_3$ of~$p$ in $\R$ 
and an event $\Omega_2$ with full probability
such that for every sequence $k_n$
with $k_n/n \to \rho \in U_3$ and all $\h \in \Omega_2$, 
$$
Q_n(k_n) = \sqrt{\frac{-\al''(\rho)}{2\pi n }}\;
           W_\8(\be(\rho)) \; \exp\Big\{n \al\Big(\frac{k_n}{n}\Big)\Big\}\; (1+o(1))
$$
where $o(1)$ tends to 0 as $n \to \8$, and $\be(\rho)= \ln \frac{(1-p)\rho}
{p(1-\rho)}$. The neighborhood~$U_3$ is contained in~$\V(p)$, hence we have $\al=-\lh^*$ with
$\lh^*$ given by (\ref{eq:lh*}).
\end{theo}
We note that the leading order is deterministic, but the prefactor 
is random (as $ W_\8$), depending on the particular realization 
of the Bernoulli field. This theorem is a corollary of a more refined result
(Theorem~\ref{th:sharp2}), which can be found in Section~\ref{sec:sharp}.
This will be proved by complex analysis arguments, considering the
 Fourier transform of $H_n$ under some (polymer) measure. 
 Fourier methods are quite strong, they are used in a different
spirit in \cite{BMP} to obtain sharp results on the polymer path itself
for small $\be$. The disadvantage is that we have to restrict the parameter 
domain. It would be tempting to use only real variable techniques 
as in the Ornstein-Zernike theory for the Bernoulli bond percolation 
\cite{CI}, but we take another, shorter route.

\goodbreak

\begin{rem}\end{rem} \vspace{-3mm}
The model is also interesting with real-valued $\h(t,x)$ with
general distribution. This is motivated by first-passage time percolation.
Our results at the exponential order remain
valid for variables with exponential moments. For the case of the
Gaussian law, we mention the recent preprint \cite{Ku} on the
so-called REM conjecture: it is proved that the local statistics of 
$(H_n(S); S \in \calP_n)$ approach that of a Poisson point process,
provided that one focuses on values distant from the mean $\E H_n$ 
by at most $o(n^{1-\eps})$.
%%%%%%% conjecture: in discrete case, does $Q_n(k_n)/\EQ_n(k_n)$ converge?

We can interpret our last result in this spirit. In our case, 
$(H_n(S); S \in \calP_n)$ spreads on the lattice, and natural local statistics
of the energy levels are the ratios $Q_n(k_n)/\E Q_n(k_n)$. For $d \geq 3$
and $k_n \sim n \rho \in U_3$, 
$$Q_n(k_n)/\E Q_n(k_n) \simeq W_\8(\be(\rho)) 
$$
since $\E W_\8(\be)=1$.
We emphasize that here the energy level~$k_n$ is of order~$n$
(far from the bulk), and that the 
limit is not universal but depends on the lattice and the law of the 
environment $\h$. 

\section{Logarithmic asymptotics} \label{sec:logas}

\subsection{Proof of Theorem~\ref{th:Rn}} \label{sec:pfth:Rn}

We start by introducing 
%% the polymer 
some probability measures.

Let $P$ be the law of the 
simple random walk on 
$\Z^d$ starting from 0, i.e. the probability measure on the space $\calP_\8$
of infinite paths making the increments $S_t-S_{t-1}$ independent and 
uniformly distributed on the set of $2d$ neighbors of $0 \in \Z^d$. 
Observe
that the restriction of $P$ to  paths of length $n$ is the normalized
counting measure on $\calP_n$, and so the partition function takes 
now the familiar form ($E_P$ is the expectation with respect to~$P$)
\begin{equation} \label{eq:ZnP}
Z_n = (2d)^n E_P [  \exp\{ \be H_n(S)\}].
\end{equation}
%    The  model of directed polymers consists in a 
%    sequence of random probability measures on $\calP_n$: the 
%    {\bf polymer measure} $\mu_n=\mu_n^{\be, \h}$ 
%    at temperature inverse $\be \in \R$
%    in the environment $\h$ is defined by
%    \begin{equation} \label{def:mu}
%    \mu_n (S)= \frac{ 1}{  Z_n} \exp\{ \be H_n(S)\} \;,\quad S \in 
%    \calP_n
%    \end{equation}
The law $\nu_n=\nu_n^{ \h}$ of $(1/n)H_n$ under $P$, given by $\nu_n(\{\rho\}):=
P(H_n(S)=n\rho)$, is such that
%% $\nu_n(\rho)=0$ if $n\rho \notin \{0,1\ldots n\}$, and 
\begin{equation} \label{eq:nun=Qn}
\nu_n(\{\rho\})= \frac{Q_n(n\rho)}{(2d)^n}
\quad {\rm if} \quad n\rho \in \{0,1,\ldots, n\}.
\end{equation}
We extend  $\nu_n$ to a probability  measure on $\R$ that we still denote
by $\nu_n$,
%% (with a slight abuse of notation),
$ \nu_n(A) = \sum_{\rho \in A, n\rho \in \{0,1,\ldots, n\}} \nu_n(\{\rho\}),
A \subset \R$.

All what we need to obtain the proof of Theorem~\ref{th:Rn}, is to prove
that $\nu_n$ obeys an almost sure large deviation 
principle, see Proposition~\ref{LDPnu} below. Recall first the event $\Omega_0(\be)$ 
from~(\ref{def:phips}), and define the event $\Omega_0 = \bigcap_{\be \in \Q} 
\Omega_0(\be)$. Then, we have $\IP( \Omega_0 )=1$, and on this event 
the convergence (\ref{def:phips}) holds for any real number~$\be$ by convexity and monotonicity.

\begin{prop} \label{LDPnu}
The function  
$$
I(\rho)=\ln(2d)+\phi^*(\rho)
 \in [0, \ln(2d)] \cup \{+\8\}
$$ 
is lower semi-continuous and convex
on $[0,1]$.
Moreover, for all $\h \in \Omega_0$ the sequence 
$(\nu_n, n\geq 1)$ obeys a large deviation 
principle with rate function $I$. That is,
\begin{itemize}
\item[(i)] for any closed $F \subset [0,1]$, we have 
$$
\limsup_{n \to \8}n^{-1}\ln \nu_n(F) \leq - \inf_{\rho \in F} I(\rho),
$$
\item[(ii)] for any open (in the induced topology on $[0,1]$) $G \subset [0,1]$, we have  
$$
\liminf_{n \to \8}n^{-1}\ln \nu_n(G) \geq - \inf_{\rho \in G} I(\rho).
$$
\end{itemize}
\end{prop}

Now, we first finish the proof of Theorem~\ref{th:Rn}, and then prove the above proposition.

\noindent
{\it Proof of Theorem \ref{th:Rn}}.
Assume that $\rho \in [p,\rho^+)$ and $\h \in \Omega_0$.
Applying (i) of Proposition \ref{LDPnu} with $F=[\rho, 1]$
and using (\ref{eq:nun=Qn}) together with
the fact that  $\rho \geq p$,
we see that the limit in (\ref{th:eqRn})
is not larger than $\ln(2d)-I(\rho)=\al(\rho)$. 
Applying (ii) of Proposition \ref{LDPnu} with $G=(\rho+\eps, 1]$ ($\eps>0$)
and using the fact that  $\rho \geq p$,
we see that the limit is at least $\al(\rho+\eps)$. Since $\rho < \rho^+$,
 this quantity tends to $\al(\rho)$ as
$\eps \searrow 0$. This proves (\ref{th:eqRn}) for $\rho \in [p, \rho^+)$.
The case  $\rho \in (\rho^-,p)$ is completely similar. 
% For $\rho \geq \rho^+$,
% we see obtain from (i) of Proposition \ref{LDPnu} that a.s., for all $\eps>0$
% it holds $\nu_n([\rho, 1])\leq e^{-n[I(\rho)-\eps]}$ for large enough $n$.
% This implies the desired result for $\rho=\rho^+$. 
Finally, when $\rho>\rho^+$ (the case  $\rho < \rho^-$ is similar)
we have $I(\rho)=\8$ and then $R_n(\rho)=0$ for large $n$, proving~(\ref{th:eqRn}) in this case. 
\qed

\noindent
{\it Proof of Proposition \ref{LDPnu}}.
The properties of $I$ are clear from the definition. 

 Fix $\h \in \Omega_0$. In view of (\ref{def:phips})  and (\ref{eq:ZnP}),
 the Laplace transforms of $\nu_n(\cdot)=P((1/n)H_n=\cdot)$ 
have logarithmic asymptotics:
$$
\lim_{n \to \8}   \frac{1}{n} \ln E_P (\exp \{\be H_n(S)\}) = 
\phi(\be) -\ln(2d)
$$
for all real $\be$. From the G\"artner-Ellis theorem (Theorem~2.3.6 in~\cite{DZ}), 
the full statement~(i) in Proposition~\ref{LDPnu} follows, and we obtain
for open $G \subset [0,1]$ that 
\begin{equation} \label{eq:eindhoven}
\liminf_{n \to \8}  \frac{1}{n}
\ln \nu_n(G) \geq - \inf\{ I(\rho);
\rho \in G \cap \calE \}\;,
\end{equation}
%%it follows statement (i) of Theorem \ref{LDPnu} holds, though  statement 
%% (ii) holds with the infimum taken over $\rho \in G \bigcap \calE$, 
where 
$$
\calE = \Big\{ \rho \in [0,1]: \exists \beta \;  \forall r \neq \rho, \;
\beta \rho - \phi^*(\rho) >
\beta r - \phi^*(r) 
\Big\}
$$
is the set of exposed points of $\phi^*$ from (\ref{def:phi*}).
Its complement is the set of all points $\rho$ such that $\phi^*$ is linear 
in a neighborhood of $\rho$.  We will improve~(\ref{eq:eindhoven})
into (ii) of Proposition \ref{LDPnu} with a subadditivity argument.
We start by showing that $\phi$ is differentiable at 0 with 
that $\phi'(0)=p$. Indeed, using Jensen inequality twice, we have
\begin{equation*}
\frac{1}{n} \ln E_P[e^{\E \be (H_n-np)}] 
\leq
\E \frac{1}{n} \ln E_P[e^{ \be (H_n-np)}]
\leq  \frac{1}{n} \ln \E E_P[e^{ \be (H_n-np)}]
\end{equation*}
Computing the extreme terms and taking the limit $n \to \8$ for 
the middle one,  we get
\begin{equation*}
0 \leq \phi(\be) - \be p
\leq  \la(\be)-\be p\;,
\end{equation*}
which shows that $\phi'(0)=p$ since $\la'(0)=p$. This implies that $p \in 
\calE$ and that $\calE$ is a neighborhood of $p$. 
Let $\rho \in (\rho^-, \rho^+) \cap G$ be a non-exposed point of 
$\phi^*$. For definiteness, we assume $\rho >p$. Let 
$$\rho_1=\sup \{
\rho' \in \calE; \rho'<\rho\}\;,\quad \rho_2=\inf \{
\rho' \in \calE; \rho'>\rho\}\;.
$$
Recall that $\phi$ is strictly convex by Theorem \ref{theo:phistrictcvx}
-- that we will prove below independently. This implies that 
 the function~$\phi^*$  cannot have a linear piece
that goes up to~$\rho^+$, cf.~Figure~\ref{fig:phi}.
Then, $p < \rho_1< \rho < \rho_2< \rho^+$, and $ \rho_1, \rho_2 \in \calE$.
Let $\gamma \in (0,1)$ such that $\rho=\gamma \rho_1 + (1-\gamma) \rho_2$. 
Since the interval $(\rho_1, \rho_2)$ consists of non-exposed points,
we have $I(\rho)=\gamma I(\rho_1) + (1-\gamma) I(\rho_2)$.
Since~$G$ is open and contains~$\rho$, we can find  
$\eps>0$ and $k,\ell \in \N^*$ such that 
$$
|u-\rho_1|<\eps, |v-\rho_2|<\eps
\Longrightarrow 
\frac{ku+\ell v}{k+\ell} \in G^\eps
$$
with $G^\eps$ the set of $r \in G$ at distance at least $\eps$ from the
outside of $G$.
The key fact is
  \begin{eqnarray*}
\lefteqn{\Card \Big\{ S \in \calP_{n(k+\ell)}: \frac{H_{n(k+\ell)}(S)}{n(k+\ell)} \in 
G^\eps \Big\}} \\ 
&\geq& \sum_{x \in \Z^d}
\Card \Big\{ S \in \calP_{n(k+\ell)}: \frac{H_{nk}(S)}{nk} \in 
(\rho_1-\eps, \rho_1+\eps)
, S_{nk}=x \Big\} \\ 
&& {}\times
\Card \Big\{ S \in \calP_{n(k+\ell)}: S_{nk}=x, \\
 && ~~~~~~~~~~~~~~~\frac{H_{n(k+\ell)}(S)-H_{nk}(S)}{n\ell} \in 
(\rho_2-\eps, \rho_2+\eps) \Big\}
 \\ &\geq &
\Card \Big\{ S \in \calP_{n(k+\ell)}: \frac{H_{nk}(S)}{nk} \in 
(\rho_1-\eps, \rho_1+\eps)
\Big\} \\ 
&& {} \times \min_{\|x\|_1 \leq nk}
\Card \Big\{ S \in \calP_{n(k+\ell)}: S_{nk}=x, \\
  && ~~~~~~~~~~~~~~~~~~~~~~~~~~\frac{H_{n(k+\ell)}(S)-H_{nk}(S)}{n\ell} \in 
(\rho_2-\eps, \rho_2+\eps) \Big\}
 \\ &=& 
\Card \Big\{ S \in \calP_{nk}: \frac{H_{nk}(S)}{nk} \in 
(\rho_1-\eps, \rho_1+\eps)
\Big\} \\ 
&&{}\times \min_{\|x\|_1 \leq nk}
\Card \Big\{ S \in \calP_{n\ell}: \frac{H_{n\ell}^{(nk,x)}(S)}{n\ell} \in 
(\rho_2-\eps, \rho_2+\eps) \Big\}
  \end{eqnarray*}
with $H_{n\ell}^{(nk,x)}(S)=\sum_{t=1}^n \h(t+nk,S_t+x)$ the Hamiltonian in
the time-space shifted environment. Similarly, we denote by 
$\nu_{n\ell}^{(nk,x)}$ the measure $\nu_{n\ell}^{(nk,x)}(\cdot)=
P(H_{n\ell}^{(nk,x)} \in \cdot)$.
The above display implies that
  \begin{eqnarray*}
\lefteqn{\liminf_{n \to \8}
\frac{1}{n(k+\ell)} \ln 
\nu_{n(k+\ell)}(G^\eps)}\\
 &\geq& \frac{k}{k+\ell}
\liminf_{n \to \8}
\frac{1}{nk} \ln 
\nu_{nk}\big((\rho_1-\eps, \rho_1+\eps)\big)\\
&& {}+ \frac{\ell}{k+\ell}
\liminf_{n \to \8}
\frac{1}{n\ell} \min_{\|x\|_1 \leq nk}
\ln 
\nu_{n\ell}^{(nk,x)}\big((\rho_2-\eps, \rho_2+\eps)\big)
  \end{eqnarray*}
It is straightforward to check that 
$$\liminf_{n \to \8}
\frac{1}{n(k+\ell)} \ln 
\nu_{n(k+\ell)}(G^\eps) \leq \liminf_{n \to \8}
\frac{1}{n} \ln 
\nu_{n}(G)\;,
$$
and it is not difficult to see that
  \begin{equation} \label{eq:hotel}
\liminf_{n \to \8}
\frac{1}{n\ell} \min_{\|x\|_1 \leq nk}
\ln 
\nu_{n\ell}^{(nk,x)}\big((\rho_2-\eps, \rho_2+\eps)\big)
\geq -I(\rho_2) \;,\quad \text{$\IP$-a.s.}
  \end{equation}
We postpone the  proof of (\ref{eq:hotel}) for the moment.
Hence, the key inequality implies
  \begin{eqnarray*}
 \liminf_{n \to \8}
\frac{1}{n} \ln 
\nu_{n}(G) 
&\geq &
- \frac{k}{k+\ell} I(\rho_1+\eps)-
\frac{\ell}{k+\ell} I(\rho_2+\eps)\;,\\
 \liminf_{n \to \8}
\frac{1}{n} \ln 
\nu_{n}(G) 
&\geq & - \big[\gamma I(\rho_1) + (1-\gamma) I(\rho_2)\big]= -I(\rho),
  \end{eqnarray*}
letting $\eps \searrow 0$ and $k/(k+\ell) \to \gamma$. This yields
 statement (ii) in Proposition~\ref{LDPnu}.

Now, let us prove~(\ref{eq:hotel}). By a standard concentration 
inequality (e.g., Theorem 4.2 in \cite{CV}), we have
$$
\IP( |\ln Z_n - \E \ln Z_n| \geq u) \leq 2 \exp \Big\{- \frac{u^2}{4 \be^2 n}\Big\}\;.
$$
Therefore we have, $\IP$-a.s. as $n \to \8$,
$$
\max_{ \|x\|_1 \leq m \leq n} \Big\vert \frac{1}{n} \ln Z_n^{(m,x)}(\be) - 
\phi(\be)
 \Big\vert \to 0 \; ,\quad \be \in \R\;,
$$
with $Z_n^{(m,x)}$ the partition function associated to
$H_{n}^{(m,x)}$. Since $\rho_2$ is an exposed point for $\phi^*$,
(\ref{eq:hotel}) follows from   the G\"artner-Ellis theorem.
\qed

Let us comment on the above proof.  We could improve (\ref{eq:eindhoven})
into the full lower bound (ii) in Proposition \ref{LDPnu} 
with a subadditivity argument, implying convexity of the rate function.
If we knew that  $ (\rho^-,
\rho^+) \subset \calE$ -- or, equivalently, that $\phi$ is differentiable --, 
we could directly conclude without this extra argument.
We tried to prove it, but we could not. We state it as a conjecture:
\begin{conj} \label{phidiff}
The function $\phi$ is everywhere differentiable. 
\end{conj}

\subsection{Proof of Theorem \ref{th:al-lh}} \label{sec:al-lh}

   \begin{enumerate}
   \item By Jensen inequality,
$$
 \E \ln Z_n \leq  n \lh(\be)\;.
$$
Then, $\phi(\be) \leq \lh(\be)$, which implies $\phi^*(\rho) \geq \lh^*(\rho)$
from the definition of Legendre transform.
The inequality now follows from $\al \leq -\phi^*$. 
   \item Set $\phi_n(\beta)= n^{-1}\E \ln Z_n (\beta)$.
From Theorem 1.1 in \cite{CY} we have
$$ \lh'(\beta) \geq  \phi_n'(\beta)  $$
for all $\beta \geq 0$. Hence, for $\rho \geq p$, the reciprocal functions 
are such that 
$$ (\lh')^{-1}(\rho) \leq  (\phi_n')^{-1}(\rho)\;. $$
Since $(\lh')^{-1}=(\lh^*)'$ and $ (\phi_n')^{-1}=(\phi_n^*)'$, we have 
$$(\lh^*)'(\rho) \leq  (\phi^*)'(\rho)$$
for all $\rho \geq p$ where $\phi^*$ is differentiable. Since $\alpha=-\phi^*$
for $\rho \neq \rho^+$, this proves 
the first half of the desired statement. The other half is similar.
   \item From Theorem 1.1 in \cite{CY} it is known that the set
$$
\W(p) = \{ \beta \in \R: \phi(\be)=\lh(\be)\}
$$
is an interval containing 0.
Let $\be \in \W(p)$, and $\rho=\la'(\be)=
\lh'(\be)$.  From Theorem 2.3~(a) in~\cite{CSY} it is known that
$\beta \in \W(p)$ implies $\phi^*(\rho)\leq 0$.
Then, the supremum defining  $\lh^*(\rho)$ is achieved at 
$\be$, which implies the first equality in
  \begin{eqnarray*}
-\lh^*(\rho) = -[ \be \rho - \lh(\be)] 
= -  [ \be \rho - \phi(\be)] 
= -\phi^*(\rho) =\al(\rho)\;,
  \end{eqnarray*}
where the second equality holds for $\be \in \W(p)$, the third one 
because of $ \phi'(\be)=\lh'(\be)=\rho$, and the last one because
$\phi^*(\rho)\leq 0$.

Let now $\be \notin \W(p)$, and $\rho=\la'(\be)$. Then,
  \begin{eqnarray*}
-\lh^*(\rho) = -[ \be \rho - \lh(\be)] 
> -  [ \be \rho - \phi(\be)] 
\geq -\phi^*(\rho) \geq \al(\rho)\;.
  \end{eqnarray*}
Observe that $\la'$ is a diffeomorphism from $\R$ to $(0,1)$. From 
this we can identify the set $\V(p)$ defined by (\ref{def:Vp}),
\begin{equation}
  \label{eq:Vp}
  \V(p)= \{ \la'(\be); \be \in \W(p)\}\;,
\end{equation}
which is an interval containing $p$.
\item When $d=1$, it is known that $\W(p)=\{0\}$, see Theorem 1.1 in \cite{CY}.
Hence, $\V(p)$ reduces to~$\{p\}$.
\item When $d\geq 3$, from celebrated results of Imbrie and Spencer \cite{IS}, 
Bolthausen \cite{Bolt},
it is known that $\W(p)$ contains a neighborhood of 0.
In view of~(\ref{eq:Vp}), $\V(p)$ is in its turn a neighborhood of $p$.
\item This is a consequence of \cite[example 2.1.1]{CSYrev}, which shows
for instance that, if $p >\pi_d$, then $\W(p)\supset \R^+$. Indeed, in view of
(\ref{eq:Vp}), this implies that $ \V(p)$ contains $[p,1)$, and
$\al$ is still equal to $-\lh^*$ at $\rho=1$ by upper semi-continuity
of both functions. The case of  $p <1-\pi_d$ is similar.
\item This is a consequence of \cite[example 2.2.1]{CSYrev}, which shows
for instance that, if $p < (1/2d)$, then $\W(p)$ is bounded from above. The other 
case is similar.
   \end{enumerate}
\qed

\subsection{Strict convexity of the free energy} 

The aim of this section is to prove Theorem~\ref{theo:phistrictcvx}.
We start with  a variance estimate  analogous
to that for Gibbs field in  \cite{DN}.

\begin{lmm} \label{lem:phistrictcvx}
For any compact set $K \subset \R$, there exists a positive
 constant $C=C_K$ such that
$$
\E (\ln Z_n)''(\be) \geq Cn \;, \qquad \be \in K
$$
\end{lmm}

\noindent
{\it Proof:} The polymer measure at inverse temperature $\be$
with environment $\h$ is the random probability measure 
$\mu_n=\mu_n^\be$ on the path space defined by 
\begin{equation}
  \label{eq:mun}
\mu_n (\{S\})= Z_n^{-1} \exp\{ \be H_n(S)\} \;, \quad
S=(S_1,\ldots, S_n) \in {\mathcal P}_n
\;.
\end{equation}
For simplicity we write $\mu_n (S_1,\ldots, S_n)$
for $\mu_n (\{(S_1,\ldots, S_n)\})$.
The polymer measure is Markovian (but time-inhomogeneous), and
$$
 (\ln Z_n)''(\be) = {\rm Var}_{\mu_n} (H_n).
$$
Let $\Si_t$ be the $t$-coordinate mapping on $\calP_n$ given by  
$\Si_t(S)=S_t$, and regard it as a random variable. 

Define 
\begin{equation}
  \label{eq:calI}
  \calI (x,y) =\{z \in \Z^d: \|x-z\|_1=\|z-y\|_1=1\}\;,\quad
x,y \in \Z^d
\end{equation}
the set of 
lattice points which are next to  both $x$ and $y$. 
The set $ \calI (x,y)$ is empty 
except if $y$ can be reached 
in two steps by the simple random walk from $x$; in this case
its cardinality is equal to $2d, 2$ or 1 according to $y=x, \|y-x\|_\8=1$
or $ \|y-x\|_\8=2$. 
The Markov property implies that, under
$\mu_{2n}$, $\Si_1,\Si_3, \ldots, \Si_{2n-1}$ are independent
conditionally on $\Si^e:=(\Si_2,\Si_4, \ldots, \Si_{2n})$,
and the  law of $\Si_{2t-1}$ given $\Si^e$ only depends
on $\Si_{2t-2}, \Si_{2t}$, and has support $\calI (\Si_{2t-2}, \Si_{2t})$.

From the variance decomposition under conditioning, we have
\begin{eqnarray*}
  {\rm Var}_{\mu_{2n}} (H_{2n}) 
&=&
  E_{\mu_{2n}} {\rm Var}_{\mu_{2n}} (H_{2n} \mid \Si^e) 
 +  {\rm Var}_{\mu_{2n}} ( E_{\mu_{2n}}[H_{2n} \mid \Si^e])\\
&\geq& 
 E_{\mu_{2n}} {\rm Var}_{\mu_{2n}} (H_{2n} \mid \Si^e) \\
&=&
 E_{\mu_{2n}} {\rm Var}_{\mu_{2n}} ( \sum_{t=1}^n
\h(2t-1, \Si_{2t-1})
\mid \Si^e) \\
&=&
\sum_{t=1}^n  E_{\mu_{2n}} {\rm Var}_{\mu_{2n}} ( 
\h(2t-1, \Si_{2t-1})
\mid \Si^e)
\end{eqnarray*}
where $ E[\;\cdot \mid \Si^e],  {\rm Var} (\;\cdot \mid \Si^e) $ denote
conditional expectation and conditional variance. To obtain the last equality
we used the conditional independence.
Define the event
$$
M(\h,t,y,z)= \Big\{
\Card \big\{ \h(t, x); x \in \calI (y,z)\big\}
= 2\Big\}\;.
$$
The reason for introducing $M(\h,t,y,z)$ is that on this event,
a path $S$ conditioned on $S_{t-1}=y, S_{t+1}=z$, has the option
to pick up a $\h(t,S_t)$ value that can be either 0 or 1, bringing therefore 
some amount of randomness. This event plays a key role here, as well as in 
the proof of Lemma \ref{lem:locality} below.
Note for further purpose that
\begin{equation}
  \label{eq:PM}
\IP\Big( M(\h,t,y,z) \Big) = 1 - \left( p^{\Card \calI (y,z)} +
 (1-p)^{\Card \calI (y,z)} \right)
=:q(y-z)\;.
\end{equation}
The key observation is, for all $t\leq n$ and $ \be \in 
K$,
\begin{equation}
  \label{eq:keyvar}
{\rm Var}_{\mu_{2n}} ( 
\h(2t-1, \Si_{2t-1})
\mid \Si^e) \geq C \1{M(\h,2t-1, \Si_{2t-2},\Si_{2t}) }\;,
\end{equation}
where the constant $C$ depends only on $K$ and the dimension $d$.
Indeed, on the event $M(\h,2t-1, S_{2t-2},S_{2t})$,  the variable
$\h(2t-1, \Si_{2t-1})$
 brings some fluctuation under the conditional law:
it takes values 0 and 1 with probability uniformly 
bounded away from 0
provided~$\be$ remains in the compact.
Hence,
\begin{eqnarray*}
\E  {\rm Var}_{\mu_{2n}} (H_{2n}) &\geq& C 
\E  \sum_{t=1}^n \mu_{2n}[M(\h,2t-1, \Si_{2t-2},\Si_{2t})]\\
 &=& C 
\E  \sum_{t=1}^n \sum_{x,y \in \Z^d} \!\mu_{2n}( \Si_{2t-2}=x,\Si_{2t}=y)
\1{M(\h,2t-1,x,y)}
\end{eqnarray*}
For $1 \leq i \leq n$, let $\tilde \mu_{n}^{(i)}$ 
be the polymer measure in the 
environment $\tilde \h(t,x)=\h(t,x)$ if $t \neq i$, 
 $\tilde \h(i,x)=0$ for all $x$. Obviously,
$$
 C^- \tilde \mu_{n}^{(i)}(S) \leq
\mu_{n}(S) \leq
C^+ \tilde \mu_{n}^{(i)}(S)\;,\quad S \in \calP_n,
$$
with positive finite $C^-, C^+$ not depending on $n, \h, \be \in K$. Then,
with $C'=CC^-$,
\begin{eqnarray*}
\lefteqn{\E  {\rm Var}_{\mu_{2n}} (H_{2n})}\\
 &\geq& C' 
\E  \sum_{t=1}^n \sum_{x,y \in \Z^d} \tilde \mu_{2n}^{(2t-1)}
( \Si_{2t-2}\!=\!x,\Si_{2t}\!=\!y)
\1{M(\h,2t-1,x,y)}\\
&=&
C' 
\E  \sum_{t=1}^n \sum_{x,y \in \Z^d} \tilde \mu_{2n}^{(2t-1)}
( \Si_{2t-2}\!=\!x,\Si_{2t}\!=\!y)
\IP( M(\h,2t-1,x,y))
%%\qquad {\rm(independance)}\\
\\
&\geq &
2C'p(1-p)
\E  \sum_{t=1}^n \sum_{x,y \in \Z^d} \tilde \mu_{2n}^{(2t-1)}
( \Si_{2t-2}\!=\!x,\Si_{2t}\!=\!y)
 \1{\|x-y\|_\8 \leq 1}
\\
&\geq &
C'p(1-p)
\E  \sum_{t=2}^{2n} \sum_{x,y \in \Z^d} \tilde \mu_{2n}^{(2t-1)}
( \Si_{t-2}\!=\!x,\Si_{t}\!=\!y)
 \1{\|x-y\|_\8 \leq 1}\;,
\end{eqnarray*}
since we can repeat the same procedure, but conditioning on the
path at odd times. Finally, with $\Delta S_t:=S_t-S_{t-1}$ and
$C''=C'C^-p(1-p)$, we have for all $\be \in K, \eps >0$,
\begin{eqnarray} \nonumber
  \E  {\rm Var}_{\mu_{2n}} (H_{2n}) &\geq& C'' 
\E  E_{\mu_{2n}}
\sum_{t=2}^{2n} 
 \1{\Delta \Si_t \neq \Delta \Si_{t-1}  }\\  \label{eq:1106}
 &\geq& n C'' \eps \times \E  {\mu_{2n}}(A_{n, \eps})\;,
\end{eqnarray}
where
$$
A_{n, \eps} = \Big\{  S \in {\calP}_n: 
\sum_{t=2}^{2n} 
 \1{\Delta S_t \neq \Delta S_{t-1}  } \geq n \eps \Big\}\;.
$$
It is easy to see that the complement 
$$A_{n, \eps}^c = \Big\{\sum_{t=2}^{2n} 
 \1{\Delta S_t = \Delta S_{t-1}  } > n (2-\eps) \Big\}
$$
 of this set has
cardinality smaller than $\exp \{n \delta(\eps)\}$, with
$\delta(\eps) \searrow 0$ as $ \eps \searrow 0$. 
We bound
$$
\IP( \max\{ H_{2n}(S); S \in A_{n, \eps}^c\} \geq 2n\rho ) \leq
e^{n \delta(\eps)}  \times {\rm Prob}
({\mathcal B}(2n,p)\geq 2n\rho)\;,
$$ 
with ${\mathcal B}(2n,p)$ a binomial random variable. It follows that
there exists some $\rho(\eps)$ with $\rho(\eps) \searrow p$ as 
$\eps \searrow 0$ such that
the left-hand side is less than $\exp\{-n \delta(\eps)^{1/2}\}$.
For all $\h$ such that $ \max\{ H_{2n}(S); S \in A_{n, \eps}^c\} \leq 
2n\rho(\eps)$,  
we have the estimate
\begin{eqnarray*}
{\mu_{2n}}(A_{n, \eps}^c) 
%%= \frac{ \sum_{S \in A_{n, \eps}^c}e^{\be  H_{2n}} }{Z_{2n}}
&\leq & 
\exp\big\{2n [\be \rho(\eps) - \phi(\be) + \delta(\eps) + o(1)]\big\}\\
&\leq & \exp\big\{2n [ \phi^*(\rho(\eps)) + \delta(\eps) + o(1)]\big\}
\end{eqnarray*}
with $o(1) \to$ as $n \to \8$.  But, as $\eps \searrow 0$, 
$$\phi^*(\rho(\eps)) + \delta(\eps) \to \phi^*(p)=-\ln(2d)<0.$$ 
By continuity we can choose $\eps >0$ such that 
$\phi^*(\rho(\eps)) + \delta(\eps) \leq (-1/2)\ln(2d)$,
and $\E  {\mu_{2n}}(A_{n, \eps}) \to 1$ as $n \to \8$.
Finally, from (\ref{eq:1106}) we obtain the desired result for even~$n$.
The same computations apply to $\mu_{2n+1}$, yielding a similar bound.
This concludes the proof of Lemma~\ref{lem:phistrictcvx}.
\qed

\noindent
{\it Proof of Theorem \ref{theo:phistrictcvx}:} 
It follows from Lemma~\ref{lem:phistrictcvx} that, for $\be, \be' \in K$, 
$$
\phi(\be') \geq \phi(\be) +  (\be'-\be) \phi_r'(\be) + 
\frac{C_K}{2}(\be'-\be)^2
\;,\quad \be \leq \be'\;,
$$
with $\phi_r'$ the right-derivative, and a similar statement for 
$\be' \leq \be$.
Indeed, this inequality holds for $(1/n) \E \ln Z_n$ instead of $\phi$,
and we can pass to the limit $n \to \8$. This yields the strict convexity
of $\phi$. By a classical property of Legendre duality, it implies 
the differentiability of  $\phi^*$.   
\qed

%%%%%%%%%%%%%%%%%%%%%%%%
\section{Sharp asymptotics}
\label{sec:sharp}
%%%%%%%%%%%%%%%%

Assume $d \geq 3$.
Let $U_0$ be the open set in the complex plane given by
$
U_0=\{ \be \in \C : |\Im \be| < \pi\}.
$
Then, $U_0$ is a neighborhood of the real axis, and
$ \la (\be)= \log \E[ \exp \{\be \h(t,x)\} ]$ is an analytic function on $U_0$.
Define, for $n \geq 0$ and $ \be \in U_0$,
\begin{equation}
  \label{defmartc}
  W_n(\be) = E_P\Big[ \exp \Big( \be \sum_{t=1}^{n} \h(t,S_t)
    - n \la(\be)\Big) \Big] \;.
\end{equation}
Then, for all $ \be \in U_0$, the sequence $(W_n(\be), n \geq 0)$
is a $(\G_n)_n$-martingale with complex values, where $\G_n=\sigma\{\h(t,x); 
t \leq n, x \in \Z^d\}$. At the same time, for each $n$ and $\h$,
$W_n(\be)$ is an analytic function of $\be \in U_0$. 

Define the real subset
\begin{equation}
  \label{eq:L2}
U_1=
\big\{\be \in \R \;:\;  \la(2\be)-2\la(\be) < - \ln \pi_d
\big\}\;,
\end{equation}
which is an open interval $(\be_1^-,\be_1^+)$
containing 0 ($-\8\leq 
\be_1^-<0<\be_1^+\leq +\8$). The following is established in \cite{CY}:

\begin{prop} \label{prop:martc}
Define $U_2$ to be the connected component of the set
$$
\Big\{\be \in U_0 \;:\;   \la(2 \Re \be)-2\Re\la(
\be) < - \ln \pi_d
%%P[ \exists n > 0: S_n=0]
\Big\}
$$
which contains the origin. Then, $U_2$
is a complex neighborhood of $U_1$. Furthermore, there exists an event
$\Omega_1$ with $\IP(\Omega_1)=1$ 
such that, 
$$
 W_n(\be) \to W_\8(\be)\;{\rm as\ } n \to \infty,\quad \text{for all } \h \in \Omega_1,
\be \in U_2\;,
$$
where the convergence is locally uniform. 
In particular, the limit
 $ W_\8(\be)$ is holomorphic in $U_2$, and all derivatives of $W_n$  
converge locally 
uniformly to the corresponding ones of $ W_\8$. Finally, 
$W_\8(\be)>0$ for all $\be \in U_1$, $\IP$-a.s.
\end{prop}
%%%%%%%%%%%%%%%%%%%%%%%%%%%%%%%%%%%%%%%%%%%%%%%%%%%%%%%%%%%%
For the sake of completeness we repeat the proof here.
\medskip

\noindent
{\it Proof of Proposition \ref{prop:martc}}: Since
$\overline{(e^z)}=e^{\overline{z}}$ and $\overline{\E[f]}=
\E[\overline{f}]$, we have
$\overline{\la (\be)}=\la (\overline{\be})$, and
\begin{eqnarray}
\E\Big[ |W_n(\be)|^2\Big]&=& \E \Big[  E_P[
\exp\{ \be H_n(S) -  n  \la(\be)\}]
 E_P[\exp\{ \overline{\be} H_n(\tilde S)-  n  \overline{\la(\be)}\}]
\Big] \nonumber
\\\nonumber
&=&   E_{P^{\otimes 2}}\Big[ \E \Big[
\exp\{ \be H_n(S) + \overline{\be} H_n(\tilde S)- 2 n \Re \la(\be)\}
\Big]\Big]\\\nonumber
&=&  E_{P^{\otimes 2}}\Big[ \exp\Big\{  [\la(2 \Re \be)-2\Re\la(
\be)]
\sum_{t=1}^n \1{S_t=\tilde S_t}
\Big\}\Big]\\
& \leq & E_{P^{\otimes 2}}\Big[ \exp\Big\{ [\la(2 \Re \be)-2\Re\la(
\be)]
\sum_{t=1}^\8
\1{S_t=\tilde S_t}
\Big\}\Big] \\
&<& \8
\label{matmata}
\end{eqnarray}
if $\be \in U_2$. Indeed, the random variable $\sum_{t=1}^\8
\1{S_t=\tilde S_t}$ (which is the number of meetings between two independent 
$d$-dimensional simple 
random walks) is geometrically distributed with parameter $\pi_d$.

For any real $\be  \in U_2$, the positive martingale $W_n(\be)$ is 
bounded in $L^2$, hence it converges almost surely and in ${\mathcal L}^2$-norm
to a non-negative limit
$W_\8(\be)$. Moreover, the event $\{W_\8(\be)=0\}$ is a tail event, so it
has probability 0 or 1. Since $\E W_\8(\be)=1$, we have necessarily
$W_\8(\be)>0$, $\IP$-a.s.

We need a stronger convergence result.
Fix a point $\be \in U_2$ and a radius $r>0$ such that the closed disk
$D(\be, r) \subset U_2$. Choosing $R >r$ such that $D(\be, R)
\subset U_2$, we obtain by Cauchy's integral formula for
all $\be' \in D(\be, r)$,
$$
W_n(\be')=  \frac{1}{2i\pi}
\int_{\partial D(\be, R)} \frac{W_n(z)}{z-\be'} dz
  =  \int_0^1
\frac{W_n(\be+R e^{2i\pi u})R e^{2i\pi u}}{(\be+R e^{2i\pi u})-\be' }
du\;,
$$
hence
$$
X_n:=\sup\{ |W_n(\be')|; \be' \in  D(\be, r)\} \leq
R  \int_0^1  \frac{|W_n(\be+ R  e^{2i\pi u})|}{R-r} du\;.
$$
Letting $C=(R /(R-r))^2$, we obtain by the Schwarz inequality
  \begin{eqnarray*}
(\E[ X_n])^2 &\leq& C \E[ \int_0^1 |W_n(\be+ R  e^{2i\pi u})|^2 du]\\
    &\leq & C \sup\{  \E[ |W_n(\be'')|^2]; n \ge 1, \be'' \in  D(\be, R)\}\\
    & <& \8
  \end{eqnarray*}
in view of (\ref{matmata}). Notice now that $X_n$, a supremum of
positive submartingales,
 is itself  a positive submartingale. Since $\sup \E[ X_n]< \8$,
$X_n$ converges $\IP$-a.s. to a finite limit $X_\8$. Finally,
$$
\sup\{ |W_n(\be')|; \be' \in  D(\be, r), n \geq 1\}<\8 \quad \text{$\IP$-a.s.},
$$
 and
$W_n$ is uniformly bounded on compact subsets of $U_2$
on a set of environments of full probability. On this set,
$(W_n, n \geq 0)$ is a normal sequence \cite{rudin}
which has a unique limit on the real axis: since $U_2$ is
connected,
the full sequence converges to some limit
$W_\8$, which is holomorphic on $U_2$, and, as mentioned above, 
positive on the real axis.
 \qed

We do not know that $W_\8(\be) \neq 0$ for general
$\be \in U_2$, only for $\be \in U_1$. Therefore, for all
$\h \in \Omega_1$, we fix another complex 
neighborhood $U_3$ of $U_1$, included in
$ U_2$ and depending on $\h$, such that $W_\8$ and $W_n$ (for $n$ large)
belongs to $\C \setminus \R_-$. Recall that
\begin{equation}
  \label{eq:relZW}
Z_n(\beta)= W_n(\be) \exp\{n \lh(\be)\}
\end{equation}
by definition. 

It is sometimes convenient to consider, for  real $\be$, the $\be$-tilted
law 
$$
\nu_{n,\be}(k)= Z_n(\be)^{-1} e^{\be k} Q_n(k)\;,
\qquad k \in \{0,1,\ldots, n\}\;,
$$
which is a probability measure on the integers $0,1,\ldots,n$. Its mean is
equal to $(d/d\be) \ln Z_n(\be)$, and its variance is
\begin{equation}
  \label{eq:Dn}
D_{n,\be}= \frac{d^2}{d\be^2} \ln Z_n(\be) \;.
\end{equation}
These quantities depend also on $\h$, and $D_{n,\be}>0$ as soon as the
Bernoulli configuration ($\h(t,x), t\leq n, \|x\|_1 \leq n, 
\|x\|_1 = n\mod 2$) is not identically 0 or 1 on each ``hyperplane''
$t=k$, $k=1,\ldots, n$.
This happens eventually with probability 1, so we will not worry about
degeneracy of the variance $D_{n,\be}$. 
By positivity of the variance, for all $u$ in the range of 
 $(d/d\be) \ln Z_n(\cdot)$ there exists unique $\be=\be_n(u)\in \R$
such that
\begin{equation}
  \label{eq:ben}
 \frac{d}{d\be} \ln Z_n(\be_n(u))=u \;.
\end{equation}
Observe that the function $\be_n$ is itself random.
Define  for $\be \in \R, k \in \N$,
\begin{equation}
  \label{eq:In}
I_n(k)= \sup \{\be k - \ln Z_n(\be); \be \in \R\} - n \ln(2d)\;.
\end{equation}
(We will see in the proof of Theorem \ref{theo:sharp1} below,
that $I_n(k) \sim n I(k/n)$ with $I$ as in Proposition
\ref{LDPnu}.) For $k$ in the range of  $(d/d\be) \ln Z_n(\cdot)$,
we have
\begin{equation}
  \label{eq:In=}
I_n(k)= \be_n(k)k-\ln Z_n(\be_n(k)) - n \ln(2d)\;.
\end{equation}
Recall $(\be_1^-,\be_1^+)$ defined in (\ref{eq:L2}).
\begin{theo} \label{th:sharp2} There exist an event $\Omega_2$ with
$\IP(\Omega_2)=1$ and a real neighborhood $U_4$ of 0, $U_4 \subset
(\be_1^-,\be_1^+)$, with the following 
property.  
Let $k_n \in \{0,1,\ldots, n\}$ be a sequence such that 
$\be_n(k_n)$ remains in a compact subset~$K$ of
%% $(\be_1^-,\be_1^+)$,
 $U_4$, 
and let $\hat D_n=D_{n,\be_n(k_n)}$. Then, for all $\h \in \Omega_2$,
$$
Q_n(k_n)= \frac{1}{\sqrt{2\pi \hat D_n}} \exp \{ - I_n(k_n)+ n \ln(2d)\}
\times \big(1+o(1)\big),
$$
where $o(1) \to 0$ as $n \to \8$.
\end{theo}
\noindent
{\it Proof of Theorem \ref{th:sharp2}}.
Suppose that~$\be$ is a real number. Note that the Fourier transform of the tilted measure 
is
$$
\sum_{k=0}^n e^{iku} \nu_{n,\be}(k) = \frac{Z_n(\be+iu)}{Z_n(\be)}\;.
$$ 
    From the usual inversion formula for Fourier series we have
$$
Q_n(k_n) = Z_n(\be) e^{-\be k_n} \times \frac{1}{2\pi}
\int_{-\pi}^{\pi}   \frac{Z_n(\be+iu)}{Z_n(\be)} e^{-ik_nu} \,du\;.
$$
Taking $\be=\be_n(k_n)$ and using (\ref{eq:In=}) this becomes
\begin{equation}
  \label{eq:Fourier}
Q_n(k_n) = e^{ - I_n(k_n)+ n \ln(2d)} \times \frac{1}{2\pi}
\int_{-\pi}^{\pi}   \frac{Z_n\Big(\be_n(k_n)+iu\Big)}
{Z_n\Big(\be_n(k_n)\Big)} e^{-ik_nu} \,du\;.
\end{equation}
For the moment, $K$ is any compact subset of $(\be_1^-,\be_1^+)$. From 
the Taylor expansion of $Z_n$ at $\be=\be_n(k_n)$ and (\ref{eq:ben}), 
we have
$$
\log {Z_n\Big(\be_n(k_n)+iu\Big)}
 = \log {Z_n\Big(\be_n(k_n)\Big)} + iu k_n - \frac{u^2}{2} \hat D_n +
{\rm Rest}_n\;,
$$
where the remainder can be estimated by the Cauchy integral formula,
$$
|{\rm Rest}_n| \leq |u|^3 \delta_K^{-3} \max\{|\log Z_n(\be')|;
\be' \in D(\be'', \delta_K), \be'' \in K\}
$$
for all $|u| \leq \delta_K$, with $\delta_K>0$ equal to half of the 
distance  from $K$
to the complement of $U_3$. From Proposition \ref{prop:martc}
and the definition of $U_3$,
the above maximum is less that
$C_Kn$ for all $n\geq 1$, with $C_K$ random but finite and independent of $n$.

Moreover,  in view of Proposition \ref{prop:martc} and 
(\ref{eq:relZW},\ref{eq:Dn}), we see that
\begin{equation}
  \label{eq:equiDhn}
  \hat D_n = n \la''(\be_n(k_n))+W_n''(\be_n(k_n))
\end{equation}
is such that $ C_K'n   \leq  \hat D_n \leq C_K''n$ for some positive constants $C_K', C_K''$.

We split the integral in (\ref{eq:Fourier})
according to $|u|\leq \eps_n:=(\ln n/n)^{1/2}$ or not, and the 
first contribution is
\begin{eqnarray} \nonumber
  \lefteqn{\int_{|u|\leq \eps_n}   \frac{Z_n\Big(\be_n(k_n)+iu\Big)}
{Z_n\Big(\be_n(k_n)\Big)} e^{-ik_nu}\, du}\\
\nonumber &=&
  \int_{|u|\leq \eps_n} \exp\Big\{
- \frac{u^2}{2} \hat D_n\Big\}du (1+o(1))\\  \nonumber
&=& \frac{1}{\sqrt{\hat D_n}}
  \int_{|u|\leq \eps_n \hat D_n^{1/2}} \exp\Big\{
- \frac{u^2}{2}\Big\}du (1+o(1))\\ \label{eq:int1}
 &=& \frac{1}{\sqrt{2\pi \hat D_n}}
 (1+o(1))
\end{eqnarray}
since $\eps_n \hat D_n^{1/2} \to \8$ by (\ref{eq:equiDhn}).

Finally, to show that the other  contribution is negligible, we
need the following fact:
\begin{lmm} \label{lem:locality}
  There exist an event $\Omega_3$ with
$\IP(\Omega_3)=1$, an integer random variable $n_0$, 
a neighborhood $U_5$ of 0 in $\R$, and $\kappa>0$
such that $n_0(\h) < \8$ for $\h \in \Omega_3$ and  
$$
\Big\vert \frac{Z_n(\be+iu)}{Z_n(\be)}\Big\vert \leq \exp \{-\kappa n u^2\}
+ \exp \{-\kappa n\}
$$ 
for  $ \h \in \Omega_3, \beta \in U_5, u \in [-\pi, \pi]$, and $n \geq 
n_0(\h)$.
\end{lmm}

With the lemma to hand, for 
$ \h, \be, u$ as above,  we bound
$$
\int_{\eps_n < u \leq \pi}   \frac{Z_n\big(\be_n(k_n)+iu\big)}
{Z_n\big(\be_n(k_n)\big)} e^{-ik_nu} \,du 
= 
o\big({\hat D_n}^{-1/2}\big)
$$
where we have used $n={\mathcal O}(\hat D_n)$ of
 (\ref{eq:equiDhn}). 
Combined with (\ref{eq:int1}) and
 (\ref{eq:Fourier}) this
estimate yields the proof of the theorem, with
$\Omega_2=\Omega_1 \cap \Omega_3$, and $U_4=U_5 \cap U_3$. \qed

We turn to the proof of Lemma~\ref{lem:locality}, which states that 
the distribution $\nu_{n,\be}$ does not concentrate on a sublattice
of $\Z$, and is not too close from such a distribution. In our 
proof we take advantage of some (conditional) independance in
the variables $\h(t,S_t)$ under  $\nu_{n,\be}$. This is reminiscent 
of a construction of \cite{DT} for central limit theorem and equivalence
of ensembles for Gibbs random fields.
\medskip

\noindent
{\it Proof of Lemma \ref{lem:locality}:} In the notations
 of the proof of 
Lemma \ref{lem:phistrictcvx},
\begin{eqnarray*}
\Big\vert \frac{Z_{2n}(\be+iu)}{Z_{2n}(\be)}\Big\vert 
&=&
\Big\vert E_{\mu_{2n}}  e^{iuH_{2n}}
  \Big\vert \\
&=&
\Big\vert E_{\mu_{2n}}   E_{\mu_{2n}}\Big[ e^{iuH_{2n}}\,
\big \vert\, \Si^e \Big]
  \Big\vert \\
&\leq & 
 E_{\mu_{2n}}  \Big\vert E_{\mu_{2n}}\Big[ e^{iuH_{2n}}
\,\big \vert\, \Si^e \Big]
  \Big\vert \\
&=&
 E_{\mu_{2n}}  \prod_{t=1}^n 
\Big\vert E_{\mu_{2n}}\Big[ e^{iu\h(2t-1,S_{2t-1})}
\,\big \vert\, \Si^e \Big]
  \Big\vert
\end{eqnarray*}
by conditional independence of $\Si_1,\Si_3, \ldots, \Si_{2n-1}$ under
$\mu_{2n}$ given $\Si^e$. Recall the notation $\calI$ from (\ref{eq:calI}) and
denote by
$$
m_\ell
=
%\Card \Big\{x: \|x-S_{2t-2}am\|_1=\|S_{2t}-x\|_1=1, \h(2t-1,x)=\ell \Big\}\,,
\Card \big\{x \in \calI(S_{2t-2},S_{2t}): \h(2t-1,x)=\ell \big\}\,,
 \qquad 
\ell=0,1,\ldots,
$$
the number of sites which can be reached by the walk at time $2t-1$ and
where $\h(\cdot)$ equals to 0 and 1 respectively ($m_1+ m_0 \leq 2d$). 
Then, for $m_0, m_1 \geq 1$,
\begin{eqnarray*}
\Big| E_{\mu_{2n}}\Big[ e^{iu\h(2t+1,S_{2t+1})}
\,\big \vert \,\Si^e \Big]\Big|&=&
\Big\vert \frac{m_1e^{\be+iu}+m_0}{m_1e^{\be}+m_0}\Big\vert \\
&\leq& \exp \{- Cu^2\}\;, \qquad |u| \leq \pi\;,
\end{eqnarray*}
where the constant $C$ is uniform for $\be \in K$, and $1 \leq m_0, m_1
\leq 2d$.  We obtain
$$
\Big\vert E_{\mu_{2n}}  e^{iuH_{2n}}
  \Big\vert
\leq
 E_{\mu_{2n}} \exp\Big\{-Cu^2
\sum_{t=1}^n \1{M(\h,2t-1, \Si_{2t-2},\Si_{2t})}\Big\}\;.
$$
So far our arguments do not require $\be$ to be small. From this point,
we will use  a perturbation argument. 
Since $\mu_{2n}^{(\beta)}$ is equal to 
$P$ for $\be=0$, we study  the term 
 on the right-hand side for
the simple random walk measure $P$ instead of the polymer measure
$\mu_{2n}$, and estimate the error from this change of measure. This procedure
is rather weak, we believe that the result of the lemma holds for a much larger
range of $\be$, but we 
we do not know how to control the term in the 
right-hand side in a different way.

For $\eps >0$ we split the last expectation according to
the sum being larger or smaller than $n\eps$,
\begin{eqnarray} \nonumber
\Big\vert E_{\mu_{2n}}  e^{iuH_{2n}}
  \Big\vert
&\leq&
e^{-C\eps u^2} + \mu_{2n} \Big( 
\sum_{t=1}^n \1{M(\h,2t\!-\!1, \Si_{2t-2},\Si_{2t})} \leq n \eps \Big) \\
\label{eq:controleexp}
 &\leq&
e^{-C\eps u^2} + e^{2n\be} P \Big( 
\sum_{t=1}^n \1{M(\h,2t\!-\!1, \Si_{2t-2},\Si_{2t})} \leq n 
\eps \Big)\;\,\phantom{**}
\end{eqnarray}
by the obvious inequalities $0\leq H_{2n} \leq 2n$.  For $\gamma \in (0,1]$,
note that \goodbreak
\begin{eqnarray*}
\E \exp\!\Big\{-\gamma \1{M(\h,2t-1, S_{2t-2},S_{2t})}\Big\}
&=&
e^{-\gamma} q(S_{2t-2}-S_{2t})\\
&&{}+[1-q(S_{2t-2}-S_{2t})]\;,
\end{eqnarray*}
with $q$ defined in (\ref{eq:PM}).
Then, there exists some $C_1 >0$ such that
$$
\sup_{ \substack{x : P(\Si_2=x)>0,\\ \|x\|_\8 \leq 1}}
\big( e^{-\gamma} q(x)+[1-q(x)]\big)
\leq \exp\{ - C_1 \gamma \}\;,\qquad \gamma \in (0,1].
$$
%%    when
%%    $\|S_{2t\!-\!2}-S_{2t}\|_\8 \leq 1$, $q(S_{2t\!-\!2}-S_{2t})>0$,  
%%    and in this case the quantity
%%    in the previous display is not more than $\exp - C_1 \gamma$ ($C_1 >0$). 
Hence,
\begin{eqnarray*}
  \lefteqn{\E E_P \exp\Big\{-\gamma \sum_{t=1}^n \1{M(\h,2t-1, \Si_{2t-2},\Si_{2t})}\Big\}}\\
&=&
 E_P   \exp\Big\{-C_1 \gamma \sum_{t=1}^n \1{
\|\Si_{2t-2}-\Si_{2t}\|_\8 \leq 1}\Big\}\\
&=&
 \Big(E_P   \exp\Big\{-C_1 \gamma  \1{
\|\Si_2\|_\8 \leq 1}\Big\}\Big)^n\\
&=&
 \Big(\frac{(2d-1)e^{-C_1\gamma}+1}{2d}\Big)^n\\
&\leq &
e^{-n C_2 \gamma}
\end{eqnarray*}
with $C_2>0$. Now, we choose $\eps=C_2/2, \gamma=1$, and we get
\begin{eqnarray*}
 \lefteqn{\E  P\Big( \sum_{t=1}^n \1{M(\h,2t-1, \Si_{2t-2},\Si_{2t})} \leq n \eps\Big)}\\
&\leq&e^{n \gamma \eps}
\E E_P 
\exp\Big\{-\gamma \sum_{t=1}^n \1{M(\h,2t-1, \Si_{2t-2},\Si_{2t})}\Big\}\\
&\leq& e^{- n C_2/2} \;,
\end{eqnarray*}
and then
$$
\IP \Big(  P( \sum_{t=1}^n \1{M(\h,2t\!-\!1, \Si_{2t\!-\!2},\Si_{2t})} \leq 
n \eps)
\geq e^{-nC_2/4} \Big) \leq
e^{-nC_2/4} \;.
$$
By Borel-Cantelli lemma, the set $\Omega_3$ of all environments such that 
$$ 
P\Big( \sum_{t=1}^n \1{M(\h,2t\!-\!1, \Si_{2t\!-\!2},\Si_{2t})} \leq n \eps 
\Big)
\leq e^{-nC_2/4}\qquad {\rm  eventually},
$$ 
is of full measure. We define $n_0$ as the first integer (if exists)
from which the previous bound is fulfilled, and
$U_5=(-C_2/4,C_2/4)$. 
From 
(\ref{eq:controleexp}) we easily check 
that Lemma \ref{lem:locality} holds true with
$\kappa=\min(C\eps, C_2/2)$. \qed

\noindent
{\it Proof of Theorem \ref{theo:sharp1}:}
The theorem is a corollary of Theorem \ref{th:sharp2}, where $\Omega_2$ and 
$U_3$ are introduced. In particular we know that $\al=-\h^*$ in $U_3$.
Note that $\be (\rho)$ is the maximizer in the definition
of $\la^*(\rho)$ as a Legendre transform. 
Since $k_n/n \to \rho$, we have that $\be_n(k_n) \to \be (\rho)$. By~(\ref{eq:equiDhn}),
$\hat D_n \sim n \la''(\be (\rho))$, and by Legendre duality,
$$
 (\la^*)' \circ \la' = {\rm Id}\;,
$$
and so $\la''(\be (\rho)) = 1/(\la^*)''(\rho)$.
The only quantity left to be studied is $I_n(k_n)$. Combining (\ref{eq:In},
\ref{eq:relZW}) and performing the change of variable $\be=\be(k_n/n)+v$, 
we have
\begin{eqnarray*}
I_n(k_n)
&=&
\sup\{ \be k_n - n \lh(\be) - \ln W_n(\be); \be \in \R \}\\
 &=&
\sup\Big\{ (\be(k_n/n)+v) k_n - n \lh(\be(k_n/n)+v) \\
&&~~~~~~~~~~~~~~~~~~~ - \ln W_n(\be(k_n/n)+v); 
v \in \R \Big\}\\
 &=&
\sup\Big\{ n \Big[ \lh(\be(k_n/n))- \lh(\be(k_n/n)+v)+
\lh'(\be(k_n/n))v\Big] \\
&& \qquad\qquad\qquad
 - \ln W_n(\be(k_n/n)+v); 
v \in \R \Big\}+ n \lh^*(k_n/n) \\
 &=&  n \lh^*(k_n/n)  - \ln W_n(\be(k_n/n))\\
&&{}+ \sup\Big\{ n \Big[ \lh(\be(k_n/n))- \lh(\be(k_n/n)+v)+
\lh'(\be(k_n/n))v\Big] \\
&& \qquad\qquad\qquad
 - \ln W_n(\be(k_n/n)+v) + \ln W_n(\be(k_n/n)); 
v \in \R \Big\} \\
& =&
 n \lh^*(k_n/n)  - \ln W_n(\be(k_n/n)) +o(1) \\
& =&
 n \lh^*(k_n/n)  - \ln W_n(\be(\rho)) +o(1)
\end{eqnarray*}
by strict convexity of $\lh$ and the fact that $|\ln[W_n(\be+v)/W_n(\be)]| \leq |v|$.
 \qed

\end{document}